\documentclass{amsart}
\usepackage{amsfonts}

\newtheorem{theorem}{Theorem}
\theoremstyle{plain}

\newtheorem{corollary}{Corollary}

\newtheorem{lemma}{Lemma}

\newtheorem{proposition}{Proposition}
\newtheorem{remark}{Remark}

\numberwithin{equation}{section}
\input{tcilatex}

\begin{document}
\title[Gr\"{u}ss' Type Inequalities ]{Some Gr\"{u}ss' Type Inequalities in
Inner Product Spaces}
\author{S.S. Dragomir}
\address{School of Computer Science \& Mathematics\\
Victoria University, Melbourne, Victoria, Australia}
\email{sever@matilda.vu.edu.au}
\urladdr{http://rgmia.vu.edu.au.SSDragomirWeb.html}
\date{February 21, 2003}

\begin{abstract}
Some new Gr\"{u}ss type inequalities in inner product spaces and
applications for integrals are given.
\end{abstract}

\maketitle

\section{Introduction}

In \cite{SSD}, the author has proved the following Gr\"{u}ss' type
inequality in real or complex inner product spaces.

\begin{theorem}
\label{1.1}Let $\left( H,\left\langle .,.\right\rangle \right) $ be an inner
product space over $\mathbb{K}\left( \mathbb{K=R}\text{,}\mathbb{C}\right) $
and $e\in H,\left\Vert e\right\Vert =1.$ If $\varphi ,\gamma ,\Phi ,\Gamma $
are real or complex numbers and $x,y$ are vectors in $H$ such that the
conditions 
\begin{equation}
\func{Re}\left\langle \Phi e-x,x-\varphi e\right\rangle \geq 0\text{ and }%
\func{Re}\left\langle \Gamma e-y,y-\gamma e\right\rangle \geq 0  \label{i.1}
\end{equation}
hold, then we have the inequality 
\begin{equation}
\left\vert \left\langle x,y\right\rangle -\left\langle x,e\right\rangle
\left\langle e,y\right\rangle \right\vert \leq \frac{1}{4}\left\vert \Phi
-\varphi \right\vert \cdot \left\vert \Gamma -\gamma \right\vert .
\label{i.2}
\end{equation}
The constant $\frac{1}{4}$ is best possible in the sense that it can not be
replaced by a smaller constant.
\end{theorem}

Some particular cases of interest for integrable functions with real or
complex values and the corresponding discrete versions are listed bellow.

\begin{corollary}
\label{1.2}Let $f,g:\left[ a,b\right] \rightarrow \mathbb{K}\left( \mathbb{%
K=R}\text{,}\mathbb{C}\right) $ be Lebesgue integrable and so that 
\begin{equation}
\func{Re}\left[ \left( \Phi -f\left( x\right) \right) \left( \overline{%
f\left( x\right) }-\overline{\varphi }\right) \right] \geq 0,\;\;\;\func{Re}%
\left[ \left( \Gamma -g\left( x\right) \right) \left( \overline{g\left(
x\right) }-\overline{\gamma }\right) \right] \geq 0  \label{i.3.1}
\end{equation}%
for a.e. $x\in \left[ a,b\right] ,$ where $\varphi ,\gamma ,\Phi ,\Gamma $
are real or complex numbers and $\bar{z}$ denotes the complex conjugate of $%
z.$ Then we have the inequality 
\begin{multline}
\left\vert \frac{1}{b-a}\int_{a}^{b}f\left( x\right) \overline{g\left(
x\right) }dx-\frac{1}{b-a}\int_{a}^{b}f\left( x\right) dx\cdot \frac{1}{b-a}%
\int_{a}^{b}\overline{g\left( x\right) }dx\right\vert  \label{i.4} \\
\leq \frac{1}{4}\left\vert \Phi -\varphi \right\vert \cdot \left\vert \Gamma
-\gamma \right\vert .
\end{multline}%
The constant $\frac{1}{4}$ is best possible.
\end{corollary}

The discrete case is embodied in

\begin{corollary}
\label{1.3}Let $\mathbf{x,y\in }\mathbb{K}^{n}$ and $\varphi ,\gamma ,\Phi
,\Gamma $ are real or complex numbers so that 
\begin{equation}
\func{Re}\left[ \left( \Phi -x_{i}\right) \left( \overline{x_{i}}-\overline{%
\varphi }\right) \right] \geq 0,\func{Re}\left[ \left( \Gamma -y_{i}\right)
\left( \overline{y_{i}}-\overline{\gamma }\right) \right] \geq 0  \label{i.5}
\end{equation}
for each $i\in \left\{ 1,...,n\right\} .$ Then we have the inequality 
\begin{equation}
\left| \frac{1}{n}\sum_{i=1}^{n}x_{i}\overline{y_{i}}-\frac{1}{n}%
\sum_{i=1}^{n}x_{i}\cdot \frac{1}{n}\sum_{i=1}^{n}\overline{y_{i}}\right|
\leq \frac{1}{4}\left| \Phi -\varphi \right| \cdot \left| \Gamma -\gamma
\right| .  \label{i.6}
\end{equation}
The constant $\frac{1}{4}$ is best possible.
\end{corollary}

For other applications of Theorem \ref{1.1}, see the recent paper \cite%
{SSDIG}.

In the present paper we show that the condition $\left( \ref{i.1}\right) $
may be replaced by an equivalent but simpler assumption and a new proof of
Theorem \ref{1.1} is produced. A refinement of the Gr\"{u}ss' type
inequality $\left( \ref{i.2}\right) ,$ some companions and applications for
integrals are pointed out as well.

\section{An Equivalent Assumption}

The following lemma holds.

\begin{lemma}
\label{l2.1} Let $a,x,A$ be vectors in the inner product space $\left(
H,\left\langle .,.\right\rangle \right) $ over $\mathbb{K}\left( \mathbb{K=R}%
\text{,}\mathbb{C}\right) $ with $a\neq A.$ Then 
\begin{equation*}
\func{Re}\left\langle A-x,x-a\right\rangle \geq 0
\end{equation*}
if and only if 
\begin{equation*}
\left\Vert x-\frac{a+A}{2}\right\Vert \leq \frac{1}{2}\left\Vert
A-a\right\Vert .
\end{equation*}
\end{lemma}

\begin{proof}
Define 
\begin{equation*}
I_{1}:=\func{Re}\left\langle A-x,x-a\right\rangle ,I_{2}:=\frac{1}{4}%
\left\Vert A-a\right\Vert ^{2}-\left\Vert x-\frac{a+A}{2}\right\Vert ^{2}.
\end{equation*}
A simple calculation shows that 
\begin{equation*}
I_{1}=I_{2}=\func{Re}\left[ \left\langle x,a\right\rangle +\left\langle
A,x\right\rangle \right] -\func{Re}\left\langle A,a\right\rangle -\left\Vert
x\right\Vert ^{2}
\end{equation*}
and thus, obviously, $I_{1}\geq 0$ iff $I_{2}\geq 0$ showing the required
equivalence.
\end{proof}

The following corollary is obvious

\begin{corollary}
\label{c2.2} Let $x,e\in H$ with $\left\Vert e\right\Vert =1$ and $\delta
,\Delta \in \mathbb{K}$ with $\delta \neq \Delta .$ Then 
\begin{equation*}
\func{Re}\left\langle \Delta e-x,x-\delta e\right\rangle \geq 0
\end{equation*}
iff 
\begin{equation*}
\left\Vert x-\frac{\delta +\Delta }{2}\cdot e\right\Vert \leq \frac{1}{2}%
\left\vert \Delta -\delta \right\vert .
\end{equation*}
\end{corollary}

\begin{remark}
\label{r2.3} If $H=\mathbb{C}$, then 
\begin{equation*}
\func{Re}\left[ \left( A-x\right) \left( \bar{x}-\bar{a}\right) \right] \geq
0
\end{equation*}
if and only if 
\begin{equation*}
\left\vert x-\frac{a+A}{2}\right\vert \leq \frac{1}{2}\left\vert
A-a\right\vert
\end{equation*}
where $a,x,A\in \mathbb{C}$. If $H=\mathbb{R}$, and $A>a$ then $a\leq x\leq
A $ if and only if $\left\vert x-\frac{a+A}{2}\right\vert \leq \frac{1}{2}%
\left\vert A-a\right\vert .$
\end{remark}

The following lemma also holds.

\begin{lemma}
\label{l2.2} Let $x,e\in H$ with $\left\Vert e\right\Vert =1.$ Then one has
the following representation 
\begin{equation}
0\leq \left\Vert x\right\Vert ^{2}-\left\vert \left\langle x,e\right\rangle
\right\vert ^{2}=\inf_{\lambda \in \mathbb{K}}\left\Vert x-\lambda
e\right\Vert ^{2}.  \label{a.1}
\end{equation}
\end{lemma}

\begin{proof}
Observe, for any $\lambda \in \mathbb{K}$, that 
\begin{align*}
\left\langle x-\lambda e,x-\left\langle x,e\right\rangle e\right\rangle &
=\left\| x\right\| ^{2}-\left| \left\langle x,e\right\rangle \right|
^{2}-\lambda \left[ \left\langle e,x\right\rangle -\left\langle
e,x\right\rangle \left\| e\right\| ^{2}\right] \\
& =\left\| x\right\| ^{2}-\left| \left\langle x,e\right\rangle \right| ^{2}.
\end{align*}
Using Schwarz's inequality, we have 
\begin{align*}
\left[ \left\| x\right\| ^{2}-\left| \left\langle x,e\right\rangle \right|
^{2}\right] ^{2}& =\left| \left\langle x-\lambda e,x-\left\langle
x,e\right\rangle e\right\rangle \right| ^{2} \\
& \leq \left\| x-\lambda e\right\| ^{2}\left\| x-\left\langle
x,e\right\rangle e\right\| ^{2} \\
& =\left\| x-\lambda e\right\| ^{2}\left[ \left\| x\right\| ^{2}-\left|
\left\langle x,e\right\rangle \right| ^{2}\right]
\end{align*}
giving the bound 
\begin{equation}
\left\| x\right\| ^{2}-\left| \left\langle x,e\right\rangle \right| ^{2}\leq
\left\| x-\lambda e\right\| ^{2},\lambda \in \mathbb{K}\text{.}  \label{a.2}
\end{equation}
Taking the infimum in $\left( \ref{a.2}\right) $ over $\lambda \in \mathbb{K}
$, we deduce 
\begin{equation*}
\left\| x\right\| ^{2}-\left| \left\langle x,e\right\rangle \right| ^{2}\leq
\inf_{\lambda \in \mathbb{K}}\left\| x-\lambda e\right\| ^{2}.
\end{equation*}
Since, for $\lambda _{0}=\left\langle x,e\right\rangle ,$ we get $\left\|
x-\lambda _{0}e\right\| ^{2}=\left\| x\right\| ^{2}-\left| \left\langle
x,e\right\rangle \right| ^{2},$ then the representation $\left( \ref{a.1}%
\right) $ is proved.
\end{proof}

We are able now to provide a different proof for the Gr\"{u}ss' type
inequality in inner product spaced mentioned in Introduction, than the one
from paper \cite{SSD}.

\begin{theorem}
\label{t2.1}Let $\left( H,\left\langle .,.\right\rangle \right) $ be an
inner product space over $\mathbb{K}\left( \mathbb{K=R}\text{,}\mathbb{C}%
\right) $ and $e\in H,\left\Vert e\right\Vert =1.$ If $\varphi ,\gamma ,\Phi
,\Gamma $ are real or complex numbers and $x,y$ are vectors in $H$ such that
the conditions $\left( \ref{i.1}\right) $ hold, or, equivalently, the
following assumptions 
\begin{equation}
\left\Vert x-\frac{\varphi +\Phi }{2}\cdot e\right\Vert \leq \frac{1}{2}%
\left\vert \Phi -\varphi \right\vert ,\left\Vert y-\frac{\gamma +\Gamma }{2}%
\cdot e\right\Vert \leq \frac{1}{2}\left\vert \Gamma -\gamma \right\vert
\label{a.3}
\end{equation}%
are valid. Then one has the inequality 
\begin{equation}
\left\vert \left\langle x,y\right\rangle -\left\langle x,e\right\rangle
\left\langle e,y\right\rangle \right\vert \leq \frac{1}{4}\left\vert \Phi
-\varphi \right\vert \cdot \left\vert \Gamma -\gamma \right\vert .
\label{a.4}
\end{equation}%
The constant $\frac{1}{4}$ is best possible.

\begin{proof}
It can be easily shown (see for example the proof of Theorem 1 from \cite%
{SSD}) that 
\begin{equation}
\left\vert \left\langle x,y\right\rangle -\left\langle x,e\right\rangle
\left\langle e,y\right\rangle \right\vert \leq \left[ \left\Vert
x\right\Vert ^{2}-\left\vert \left\langle x,e\right\rangle \right\vert ^{2}%
\right] ^{\frac{1}{2}}\left[ \left\Vert y\right\Vert ^{2}-\left\vert
\left\langle y,e\right\rangle \right\vert ^{2}\right] ^{\frac{1}{2}},
\label{a.5}
\end{equation}%
for any $x,y\in H$ and $e\in H,\left\Vert e\right\Vert =1.$ Using Lemma \ref%
{l2.2} and the conditions $\left( \ref{a.3}\right) $ we obviously have that 
\begin{equation*}
\left[ \left\Vert x\right\Vert ^{2}-\left\vert \left\langle x,e\right\rangle
\right\vert ^{2}\right] ^{\frac{1}{2}}=\inf_{\lambda \in \mathbb{K}%
}\left\Vert x-\lambda e\right\Vert \leq \left\Vert x-\frac{\varphi +\Phi }{2}%
\cdot e\right\Vert \leq \frac{1}{2}\left\vert \Phi -\varphi \right\vert
\end{equation*}%
and 
\begin{equation*}
\left[ \left\Vert y\right\Vert ^{2}-\left\vert \left\langle y,e\right\rangle
\right\vert ^{2}\right] ^{\frac{1}{2}}=\inf_{\lambda \in \mathbb{K}%
}\left\Vert y-\lambda e\right\Vert \leq \left\Vert y-\frac{\gamma +\Gamma }{2%
}\cdot e\right\Vert \leq \frac{1}{2}\left\vert \Gamma -\gamma \right\vert
\end{equation*}%
and by $\left( \ref{a.5}\right) $ the desired inequality $\left( \ref{a.4}%
\right) $ is obtained.

The fact that $\frac{1}{4}$ is the best possible constant, has been shown in 
\cite{SSD} and we omit the details.
\end{proof}
\end{theorem}

\section{A Refinement of Gr\"{u}ss Inequality}

The following result improving $\left( \ref{i.1}\right) $ holds

\begin{theorem}
\label{t3.1}Let $\left( H,\left\langle .,.\right\rangle \right) $ be an
inner product space over $\mathbb{K}\left( \mathbb{K=R}\text{,}\mathbb{C}%
\right) $ and $e\in H,\left\| e\right\| =1.$ If $\varphi ,\gamma ,\Phi
,\Gamma $ are real or complex numbers and $x,y$ are vectors in $H$ such that
the conditions $\left( \ref{i.1}\right) $, or, equivalently, $\left( \ref%
{a.3}\right) $ hold, then we have the inequality 
\begin{multline}
\left| \left\langle x,y\right\rangle -\left\langle x,e\right\rangle
\left\langle e,y\right\rangle \right|  \label{3.1} \\
\leq \frac{1}{4}\left| \Phi -\varphi \right| \cdot \left| \Gamma -\gamma
\right| -\left[ \func{Re}\left\langle \Phi e-x,x-\varphi e\right\rangle %
\right] ^{\frac{1}{2}}\left[ \func{Re}\left\langle \Gamma e-y,y-\gamma
e\right\rangle \right] ^{\frac{1}{2}}.
\end{multline}
\end{theorem}

\begin{proof}
As in \cite{SSD}, we have 
\begin{equation}
\left\vert \left\langle x,y\right\rangle -\left\langle x,e\right\rangle
\left\langle e,y\right\rangle \right\vert ^{2}\leq \left[ \left\Vert
x\right\Vert ^{2}-\left\vert \left\langle x,e\right\rangle \right\vert ^{2}%
\right] \left[ \left\Vert y\right\Vert ^{2}-\left\vert \left\langle
y,e\right\rangle \right\vert ^{2}\right] ,  \label{3.2}
\end{equation}%
\begin{equation}
\left\Vert x\right\Vert ^{2}-\left\vert \left\langle x,e\right\rangle
\right\vert ^{2}=\func{Re}\left[ \left( \Phi -\left\langle x,e\right\rangle
\right) \left( \overline{\left\langle x,e\right\rangle }-\overline{\varphi }%
\right) \right] -\func{Re}\left\langle \Phi e-x,x-\varphi e\right\rangle
\label{3.3}
\end{equation}%
and 
\begin{equation}
\left\Vert y\right\Vert ^{2}-\left\vert \left\langle y,e\right\rangle
\right\vert ^{2}=\func{Re}\left[ \left( \Gamma -\left\langle
y,e\right\rangle \right) \left( \overline{\left\langle y,e\right\rangle }-%
\overline{\gamma }\right) \right] -\func{Re}\left\langle \Gamma e-x,x-\gamma
e\right\rangle .  \label{3.4}
\end{equation}%
Using the elementary inequality 
\begin{equation*}
4\func{Re}\left( a\overline{b}\right) \leq \left\vert a+b\right\vert
^{2};a,b\in \mathbb{K}\left( \mathbb{K=R}\text{,}\mathbb{C}\right)
\end{equation*}%
we may state that 
\begin{equation}
\func{Re}\left[ \left( \Phi -\left\langle x,e\right\rangle \right) \left( 
\overline{\left\langle x,e\right\rangle }-\overline{\varphi }\right) \right]
\leq \frac{1}{4}\left\vert \Phi -\varphi \right\vert ^{2}  \label{3.5}
\end{equation}%
and 
\begin{equation}
\func{Re}\left[ \left( \Gamma -\left\langle y,e\right\rangle \right) \left( 
\overline{\left\langle y,e\right\rangle }-\overline{\gamma }\right) \right]
\leq \frac{1}{4}\left\vert \Gamma -\gamma \right\vert ^{2}.  \label{3.6}
\end{equation}%
Consequently, by $\left( \ref{3.2}\right) -\left( \ref{3.6}\right) $ we may
state that 
\begin{multline}
\left\vert \left\langle x,y\right\rangle -\left\langle x,e\right\rangle
\left\langle e,y\right\rangle \right\vert ^{2}  \label{3.7} \\
\leq \left[ \frac{1}{4}\left\vert \Phi -\varphi \right\vert ^{2}-\left( %
\left[ \func{Re}\left\langle \Phi e-x,x-\varphi e\right\rangle \right] ^{%
\frac{1}{2}}\right) ^{2}\right] \\
\times \left[ \frac{1}{4}\left\vert \Gamma -\gamma \right\vert ^{2}-\left( %
\left[ \func{Re}\left\langle \Gamma e-y,y-\gamma e\right\rangle \right] ^{%
\frac{1}{2}}\right) ^{2}\right] .
\end{multline}%
Finally, using the elementary inequality for positive real numbers 
\begin{equation*}
\left( m^{2}-n^{2}\right) \left( p^{2}-q^{2}\right) \leq \left( mp-nq\right)
^{2}
\end{equation*}%
we have 
\begin{multline*}
\left[ \frac{1}{4}\left\vert \Phi -\varphi \right\vert ^{2}-\left( \left[ 
\func{Re}\left\langle \Phi e-x,x-\varphi e\right\rangle \right] ^{\frac{1}{2}%
}\right) ^{2}\right] \\
\times \left[ \frac{1}{4}\left\vert \Gamma -\gamma \right\vert ^{2}-\left( %
\left[ \func{Re}\left\langle \Gamma e-y,y-\gamma e\right\rangle \right] ^{%
\frac{1}{2}}\right) ^{2}\right] \\
\leq \left( \frac{1}{4}\left\vert \Phi -\varphi \right\vert \cdot \left\vert
\Gamma -\gamma \right\vert -\left[ \func{Re}\left\langle \Phi e-x,x-\varphi
e\right\rangle \right] ^{\frac{1}{2}}\left[ \func{Re}\left\langle \Gamma
e-y,y-\gamma e\right\rangle \right] ^{\frac{1}{2}}\right) ^{2}
\end{multline*}%
giving the desired inequality $\left( \ref{3.1}\right) .$
\end{proof}

\section{Some Companion Inequalities}

The following companion of Gr\"{u}ss inequality in inner product spaces
holds.

\begin{theorem}
\label{t4.1} Let $\left( H,\left\langle .,.\right\rangle \right) $ be an
inner product space over $\mathbb{K}\left( \mathbb{K=R}\text{,}\mathbb{C}%
\right) $ and $e\in H,\left\Vert e\right\Vert =1.$ If $\gamma ,\Gamma \in 
\mathbb{K}$ and $x,y\in H$ are so that 
\begin{equation}
\func{Re}\left\langle \Gamma e-\frac{x+y}{2},\frac{x+y}{2}-\gamma
e\right\rangle \geq 0  \label{4.1}
\end{equation}
or, equivalently, 
\begin{equation}
\left\Vert \frac{x+y}{2}-\frac{\gamma +\Gamma }{2}\cdot e\right\Vert \leq 
\frac{1}{2}\left\vert \Gamma -\gamma \right\vert ,  \label{4.2}
\end{equation}
then we have the inequality 
\begin{equation}
\func{Re}\left[ \left\langle x,y\right\rangle -\left\langle x,e\right\rangle
\left\langle e,y\right\rangle \right] \leq \frac{1}{4}\left\vert \Gamma
-\gamma \right\vert ^{2}.  \label{4.2.a}
\end{equation}
The constant $\frac{1}{4}$ is best possible in the sense that it cannot be
replaced by a smaller constant.
\end{theorem}

\begin{proof}
Start with the well known inequality 
\begin{equation}
\func{Re}\left\langle z,u\right\rangle \leq \frac{1}{4}\left\| z+u\right\|
^{2};z,u\in H.  \label{4.3}
\end{equation}
Since 
\begin{equation*}
\left\langle x,y\right\rangle -\left\langle x,e\right\rangle \left\langle
e,y\right\rangle =\left\langle x-\left\langle x,e\right\rangle
e,y-\left\langle y,e\right\rangle e\right\rangle
\end{equation*}
then using $\left( \ref{4.3}\right) $ we may write 
\begin{align}
\func{Re}\left[ \left\langle x,y\right\rangle -\left\langle x,e\right\rangle
\left\langle e,y\right\rangle \right] & =\func{Re}\left[ \left\langle
x-\left\langle x,e\right\rangle e,y-\left\langle y,e\right\rangle
e\right\rangle \right]  \label{4.4} \\
& \leq \frac{1}{4}\left\| x-\left\langle x,e\right\rangle e+y-\left\langle
y,e\right\rangle e\right\| ^{2}  \notag \\
& =\left\| \frac{x+y}{2}-\left\langle \frac{x+y}{2},e\right\rangle \cdot
e\right\| ^{2}  \notag \\
& =\left\| \frac{x+y}{2}\right\| ^{2}-\left| \left\langle \frac{x+y}{2}%
,e\right\rangle \right| ^{2}.  \notag
\end{align}
If we apply Gr\"{u}ss' inequality in inner product spaces for, say, $a=b=%
\frac{x+y}{2},$ we get 
\begin{equation}
\left\| \frac{x+y}{2}\right\| ^{2}-\left| \left\langle \frac{x+y}{2}%
,e\right\rangle \right| ^{2}\leq \frac{1}{4}\left| \Gamma -\gamma \right|
^{2}.  \label{4.5}
\end{equation}
Making use of $\left( \ref{4.4}\right) $ and $\left( \ref{4.5}\right) $ we
deduce $\left( \ref{4.2.a}\right) .$

The fact that $\frac{1}{4}$ is the best possible constant in $\left( \ref%
{4.2.a}\right) $ follows by the fact that if in $\left( \ref{4.1}\right) $
we choose $x=y,$ then it becomes $\func{Re}\left\langle \Gamma e-x,x-\gamma
e\right\rangle \geq 0,$ implying $0\leq \left\Vert x\right\Vert
^{2}-\left\vert \left\langle x,e\right\rangle \right\vert ^{2}\leq \frac{1}{4%
}\left\vert \Gamma -\gamma \right\vert ^{2},$ for which, by Gr\"{u}ss'
inequality in inner product spaces, we know that the constant $\frac{1}{4}$
is best possible.
\end{proof}

The following corollary might be of interest if one wanted to evaluate the
absolute value of 
\begin{equation*}
\func{Re}\left[ \left\langle x,y\right\rangle -\left\langle x,e\right\rangle
\left\langle e,y\right\rangle \right] .
\end{equation*}

\begin{corollary}
\label{c4.2}Let $\left( H,\left\langle .,.\right\rangle \right) $ be an
inner product space over $\mathbb{K}\left( \mathbb{K=R}\text{,}\mathbb{C}%
\right) $ and $e\in H,\left\Vert e\right\Vert =1.$ If $\gamma ,\Gamma \in 
\mathbb{K}$ and $x,y\in H$ are so that 
\begin{equation}
\func{Re}\left\langle \Gamma e-\frac{x\pm y}{2},\frac{x\pm y}{2}-\gamma
e\right\rangle \geq 0  \label{4.6}
\end{equation}
or, equivalently, 
\begin{equation}
\left\Vert \frac{x\pm y}{2}-\frac{\gamma +\Gamma }{2}\cdot e\right\Vert \leq 
\frac{1}{2}\left\vert \Gamma -\gamma \right\vert ,  \label{4.7}
\end{equation}
then we have the inequality 
\begin{equation}
\left\vert \func{Re}\left[ \left\langle x,y\right\rangle -\left\langle
x,e\right\rangle \left\langle e,y\right\rangle \right] \right\vert \leq 
\frac{1}{4}\left\vert \Gamma -\gamma \right\vert ^{2}.  \label{4.8}
\end{equation}
If the inner product space $H$ is real, then $\left( \text{for }m,M\in 
\mathbb{R}\text{, }M>m\right) $ 
\begin{equation}
\left\langle Me-\frac{x\pm y}{2},\frac{x\pm y}{2}-me\right\rangle \geq 0
\label{4.9}
\end{equation}
or, equivalently, 
\begin{equation}
\left\Vert \frac{x\pm y}{2}-\frac{m+M}{2}\cdot e\right\Vert \leq \frac{1}{2}%
\left( M-m\right) ,  \label{4.10}
\end{equation}
implies 
\begin{equation}
\left\vert \left\langle x,y\right\rangle -\left\langle x,e\right\rangle
\left\langle e,y\right\rangle \right\vert \leq \frac{1}{4}\left( M-m\right)
^{2}.  \label{4.11}
\end{equation}
In both inequalities $\left( \ref{4.8}\right) $ and $\left( \ref{4.11}%
\right) ,$ the constant $\frac{1}{4}$ is best possible.
\end{corollary}

\begin{proof}
We only remark that, if 
\begin{equation*}
\func{Re}\left\langle \Gamma e-\frac{x-y}{2},\frac{x-y}{2}-\gamma
e\right\rangle \geq 0
\end{equation*}
holds, then by Theorem \ref{t4.1}, we get 
\begin{equation*}
\func{Re}\left[ -\left\langle x,y\right\rangle +\left\langle
x,e\right\rangle \left\langle e,y\right\rangle \right] \leq \frac{1}{4}%
\left\vert \Gamma -\gamma \right\vert ^{2}
\end{equation*}
showing that 
\begin{equation}
\func{Re}\left[ \left\langle x,y\right\rangle -\left\langle x,e\right\rangle
\left\langle e,y\right\rangle \right] \geq -\frac{1}{4}\left\vert \Gamma
-\gamma \right\vert ^{2}.  \label{4.12}
\end{equation}
Making use of $\left( \ref{4.2.a}\right) $ and $\left( \ref{4.12}\right) $
we deduce the desired result $\left( \ref{4.8}\right) .$
\end{proof}

Finally, we may state and proof the following dual result as well

\begin{proposition}
\label{p4.1}Let $\left( H,\left\langle .,.\right\rangle \right) $ be an
inner product space over $\mathbb{K}\left( \mathbb{K=R}\text{,}\mathbb{C}%
\right) $ and $e\in H,\left\Vert e\right\Vert =1.$ If $\varphi ,\Phi \in 
\mathbb{K}$ and $x,y\in H$ are so that 
\begin{equation}
\func{Re}\left[ \left( \Phi -\left\langle x,e\right\rangle \right) \left( 
\overline{\left\langle x,e\right\rangle }-\overline{\varphi }\right) \right]
\leq 0,  \label{4.13}
\end{equation}%
then we have the inequalities 
\begin{align}
\left\Vert x-\left\langle x,e\right\rangle e\right\Vert & \leq \left[ \func{%
Re}\left\langle x-\Phi e,x-\varphi e\right\rangle \right] ^{\frac{1}{2}}
\label{4.14} \\
& \leq \frac{\sqrt{2}}{2}\left[ \left\Vert x-\Phi e\right\Vert
^{2}+\left\Vert x-\varphi e\right\Vert ^{2}\right] ^{\frac{1}{2}}.  \notag
\end{align}
\end{proposition}

\begin{proof}
We know that the following identity holds true $($see $\left( \ref{3.3}%
\right) )$%
\begin{equation}
\left\Vert x\right\Vert ^{2}-\left\vert \left\langle x,e\right\rangle
\right\vert ^{2}=\func{Re}\left[ \left( \Phi -\left\langle x,e\right\rangle
\right) \left( \overline{\left\langle x,e\right\rangle }-\overline{\varphi }%
\right) \right] +\func{Re}\left\langle x-\Phi e,x-\varphi e\right\rangle .
\label{4.15}
\end{equation}
Using the assumption $\left( \ref{4.13}\right) $ and the fact that 
\begin{equation*}
\left\Vert x\right\Vert ^{2}-\left\vert \left\langle x,e\right\rangle
\right\vert ^{2}=\left\Vert x-\left\langle x,e\right\rangle e\right\Vert ^{2}
\end{equation*}
by $\left( \ref{4.15}\right) $ we deduce the first inequality in $\left( \ref%
{4.14}\right) .$

The second inequality in $\left( \ref{4.14}\right) $ follows by the fact
that for any $v,w\in H$ one has 
\begin{equation*}
\func{Re}\left\langle w,v\right\rangle \leq \frac{1}{2}\left( \left\|
w\right\| ^{2}+\left\| v\right\| ^{2}\right) .
\end{equation*}
The proposition is thus proved.
\end{proof}

\section{Integral Inequalities}

Let $\left( \Omega ,\Sigma ,\mu \right) $ be a measure space consisting of a
set $\Omega ,$ a $\sigma -$algebra of parts $\Sigma $ and a countably
additive and positive measure $\mu $ on $\Sigma $ with values in $\mathbb{%
R\cup }\left\{ \infty \right\} .$ Denote by $L^{2}\left( \Omega ,\mathbb{K}%
\right) $ the Hilbert space of all real or complex valued functions $f$
defined on $\Omega $ and $2-$integrable on $\Omega ,$ i.e., 
\begin{equation*}
\int_{\Omega }\left| f\left( s\right) \right| ^{2}d\mu \left( s\right)
<\infty .
\end{equation*}

The following proposition holds

\begin{proposition}
\label{p5.1} If $f,g,h\in L^{2}\left( \Omega ,\mathbb{K}\right) $ and $%
\varphi ,\Phi ,\gamma ,\Gamma \in \mathbb{K}$, are so that $\int_{\Omega
}\left\vert h\left( s\right) \right\vert ^{2}d\mu \left( s\right) =1$ and 
\begin{align}
\int_{\Omega }\func{Re}\left[ \left( \Phi h\left( s\right) -f\left( s\right)
\right) \left( \overline{f\left( s\right) }-\overline{\varphi }\overline{%
h\left( s\right) }\right) \right] d\mu \left( s\right) & \geq 0  \label{5.1}
\\
\int_{\Omega }\func{Re}\left[ \left( \Gamma h\left( s\right) -g\left(
s\right) \right) \left( \overline{g\left( s\right) }-\overline{\gamma }%
\overline{h\left( s\right) }\right) \right] d\mu \left( s\right) & \geq 0 
\notag
\end{align}%
or, equivalently 
\begin{align}
\left( \int_{\Omega }\left\vert f\left( s\right) -\frac{\Phi +\varphi }{2}%
h\left( s\right) \right\vert ^{2}d\mu \left( s\right) \right) ^{\frac{1}{2}%
}& \leq \frac{1}{2}\left\vert \Phi -\varphi \right\vert ,  \label{5.2} \\
\left( \int_{\Omega }\left\vert g\left( s\right) -\frac{\Gamma +\gamma }{2}%
h\left( s\right) \right\vert ^{2}d\mu \left( s\right) \right) ^{\frac{1}{2}%
}& \leq \frac{1}{2}\left\vert \Gamma -\gamma \right\vert ,  \notag
\end{align}%
then we have the following refinement of Gr\"{u}ss integral inequality 
\begin{multline}
\left\vert \int_{\Omega }f\left( s\right) \overline{g\left( s\right) }d\mu
\left( s\right) -\int_{\Omega }f\left( s\right) \overline{h\left( s\right) }%
d\mu \left( s\right) \int_{\Omega }h\left( s\right) \overline{g\left(
s\right) }d\mu \left( s\right) \right\vert   \label{5.3} \\
\leq \frac{1}{4}\left\vert \Phi -\varphi \right\vert \cdot \left\vert \Gamma
-\gamma \right\vert -\left[ \int_{\Omega }\func{Re}\left[ \left( \Phi
h\left( s\right) -f\left( s\right) \right) \left( \overline{f\left( s\right) 
}-\overline{\varphi }\overline{h\left( s\right) }\right) \right] d\mu \left(
s\right) \right.  \\
\times \left. \int_{\Omega }\func{Re}\left[ \left( \Gamma h\left( s\right)
-g\left( s\right) \right) \left( \overline{g\left( s\right) }-\overline{%
\gamma }\overline{h\left( s\right) }\right) \right] d\mu \left( s\right) %
\right] ^{\frac{1}{2}}.
\end{multline}%
The constant $\frac{1}{4}$ is best possible.
\end{proposition}

The proof follows by Theorem \ref{t3.1} on choosing $H=L^{2}\left( \Omega ,%
\mathbb{K}\right) $ with the inner product 
\begin{equation*}
\left\langle f,g\right\rangle :=\int_{\Omega }f\left( s\right) \overline{%
g\left( s\right) }d\mu \left( s\right) .
\end{equation*}
We omit the details.

\begin{remark}
\label{r.1} It is obvious that a sufficient condition for $\left( \ref{5.1}%
\right) $ to hold is 
\begin{equation*}
\func{Re}\left[ \left( \Phi h\left( s\right) -f\left( s\right) \right)
\left( \overline{f\left( s\right) }-\overline{\varphi }\overline{h\left(
s\right) }\right) \right] \geq 0,
\end{equation*}%
and 
\begin{equation*}
\func{Re}\left[ \left( \Gamma h\left( s\right) -g\left( s\right) \right)
\left( \overline{g\left( s\right) }-\overline{\gamma }\overline{h\left(
s\right) }\right) \right] \geq 0,
\end{equation*}%
for $\mu -$a.e.$\;s\in \Omega ,$ or equivalently, 
\begin{align*}
\left\vert f\left( s\right) -\frac{\Phi +\varphi }{2}h\left( s\right)
\right\vert & \leq \frac{1}{2}\left\vert \Phi -\varphi \right\vert
\left\vert h\left( s\right) \right\vert \text{ \ \ and} \\
\left\vert g\left( s\right) -\frac{\Gamma +\gamma }{2}h\left( s\right)
\right\vert & \leq \frac{1}{2}\left\vert \Gamma -\gamma \right\vert
\left\vert h\left( s\right) \right\vert ,
\end{align*}%
for $\mu -$a.e.$\;s\in \Omega .$
\end{remark}

The following result may be stated as well.

\begin{corollary}
\label{c5.1}If $z,Z,t,T\in \mathbb{K}$, $\rho \in L\left( \Omega ,\mathbb{R}%
\right) ,$ $\mu \left( \Omega \right) <\infty $ and $f,g\in L^{2}\left(
\Omega ,\mathbb{K}\right) $ are such that: 
\begin{align}
\func{Re}\left[ \left( Z-f\left( s\right) \right) \left( \overline{f\left(
s\right) }-\bar{z}\right) \right] & \geq 0,  \label{5.4} \\
\func{Re}\left[ \left( T-g\left( s\right) \right) \left( \overline{g\left(
s\right) }-\bar{t}\right) \right] & \geq 0\text{ \hspace{0.05in}for a.e. }%
s\in \Omega  \notag
\end{align}%
or, equivalently 
\begin{align}
\left\vert f\left( s\right) -\frac{z+Z}{2}\right\vert & \leq \frac{1}{2}%
\left\vert Z-z\right\vert ,  \label{5.5} \\
\left\vert g\left( s\right) -\frac{t+T}{2}\right\vert & \leq \frac{1}{2}%
\left\vert T-t\right\vert \text{ \hspace{0.05in}for a.e. }s\in \Omega  \notag
\end{align}%
then we have the inequality 
\begin{multline}
\left\vert \frac{1}{\mu \left( \Omega \right) }\int_{\Omega }f\left(
s\right) \overline{g\left( s\right) }d\mu \left( s\right) \right. -\left. 
\frac{1}{\mu \left( \Omega \right) }\int_{\Omega }f\left( s\right) d\mu
\left( s\right) \cdot \frac{1}{\mu \left( \Omega \right) }\int_{\Omega }%
\overline{g\left( s\right) }d\mu \left( s\right) \right\vert  \label{5.6} \\
\leq \frac{1}{4}\left\vert Z-z\right\vert \left\vert T-t\right\vert -\frac{1%
}{\mu \left( \Omega \right) }\left[ \int_{\Omega }\func{Re}\left[ \left(
Z-f\left( s\right) \right) \left( \overline{f\left( s\right) }-\bar{z}%
\right) \right] d\mu \left( s\right) \right. \\
\times \left. \int_{\Omega }\func{Re}\left[ \left( T-g\left( s\right)
\right) \left( \overline{g\left( s\right) }-\bar{t}\right) \right] d\mu
\left( s\right) \right] ^{\frac{1}{2}}.
\end{multline}
\end{corollary}

Using Theorem \ref{t4.1} we may state the following result as well.

\begin{proposition}
\label{p5.3}If $f,g,h\in L^{2}\left( \Omega ,\mathbb{K}\right) $ and $\gamma
,\Gamma \in \mathbb{K}$ are such that $\int_{\Omega }\left\vert h\left(
s\right) \right\vert ^{2}d\mu \left( s\right) =1$ and 
\begin{equation}
\int_{\Omega }\func{Re}\left\{ \left[ \Gamma h\left( s\right) -\frac{f\left(
s\right) +g\left( s\right) }{2}\right] \cdot \left[ \frac{\overline{f\left(
s\right) }+\overline{g\left( s\right) }}{2}-\bar{\gamma}\bar{h}\left(
s\right) \right] \right\} d\mu \left( s\right) \geq 0  \label{5.7}
\end{equation}%
or, equivalently, 
\begin{equation}
\left( \int_{\Omega }\left\vert \frac{f\left( s\right) +g\left( s\right) }{2}%
-\frac{\gamma +\Gamma }{2}h\left( s\right) \right\vert ^{2}d\mu \left(
s\right) \right) ^{\frac{1}{2}}\leq \frac{1}{2}\left\vert \Gamma -\gamma
\right\vert ,  \label{5.8}
\end{equation}%
then we have the inequality 
\begin{align}
I& :=\int_{\Omega }\func{Re}\left[ f\left( s\right) \overline{g\left(
s\right) }\right] d\mu \left( s\right)  \label{5.9} \\
& \qquad \qquad \qquad -\func{Re}\left[ \int_{\Omega }f\left( s\right) 
\overline{h\left( s\right) }d\mu \left( s\right) \cdot \int_{\Omega }h\left(
s\right) \overline{g\left( s\right) }d\mu \left( s\right) \right]  \notag \\
& \leq \frac{1}{4}\left\vert \Gamma -\gamma \right\vert ^{2}.  \notag
\end{align}%
If (\ref{5.7}) and (\ref{5.8}) hold with \textquotedblleft\ $\pm $
\textquotedblright\ instead of \textquotedblleft\ $+$ \textquotedblright ,
then 
\begin{equation}
\left\vert I\right\vert \leq \frac{1}{4}\left\vert \Gamma -\gamma
\right\vert ^{2}.  \label{5.10}
\end{equation}
\end{proposition}

\begin{remark}
It is obvious that a sufficient condition for (\ref{5.7}) to hold is 
\begin{equation}
\func{Re}\left\{ \left[ \Gamma h\left( s\right) -\frac{f\left( s\right)
+g\left( s\right) }{2}\right] \cdot \left[ \frac{\overline{f\left( s\right) }%
+\overline{g\left( s\right) }}{2}-\bar{\gamma}\bar{h}\left( s\right) \right]
\right\} \geq 0  \label{5.11}
\end{equation}
for a.e. $s\in \Omega ,$ or equivalently 
\begin{equation}
\left| \frac{f\left( s\right) +g\left( s\right) }{2}-\frac{\gamma +\Gamma }{2%
}h\left( s\right) \right| \leq \frac{1}{2}\left| \Gamma -\gamma \right|
\left| h\left( s\right) \right| \text{ \hspace{0.05in}for a.e. }s\in \Omega .
\label{5.12}
\end{equation}
\end{remark}

Finally, the following corollary holds.

\begin{corollary}
\label{c5.4}If $Z,z\in \mathbb{K}$, $\mu \left( \Omega \right) <\infty $ and 
$f,g\in L^{2}\left( \Omega ,\mathbb{K}\right) $ are such that 
\begin{equation}
\func{Re}\left[ \left( Z-\frac{f\left( s\right) +g\left( s\right) }{2}%
\right) \left( \frac{\overline{f\left( s\right) }+\overline{g\left( s\right) 
}}{2}-\overline{z}\right) \right] \geq 0\text{ \hspace{0.05in}for a.e. }s\in
\Omega   \label{5.13}
\end{equation}%
or, equivalently 
\begin{equation}
\left\vert \frac{f\left( s\right) +g\left( s\right) }{2}-\frac{z+Z}{2}%
\right\vert \leq \frac{1}{2}\left\vert Z-z\right\vert \text{ \hspace{0.05in}%
for a.e. }s\in \Omega ,  \label{5.14}
\end{equation}%
then we have the inequality 
\begin{align*}
J& :=\frac{1}{\mu \left( \Omega \right) }\int_{\Omega }\func{Re}\left[
f\left( s\right) \overline{g\left( s\right) }\right] d\mu \left( s\right)  \\
& \qquad \qquad \qquad -\func{Re}\left[ \frac{1}{\mu \left( \Omega \right) }%
\int_{\Omega }f\left( s\right) d\mu \left( s\right) \cdot \frac{1}{\mu
\left( \Omega \right) }\int_{\Omega }\overline{g\left( s\right) }d\mu \left(
s\right) \right]  \\
& \leq \frac{1}{4}\left\vert Z-z\right\vert ^{2}.
\end{align*}%
If (\ref{5.13}) and (\ref{5.14}) hold with \textquotedblleft\ $\pm $
\textquotedblright\ instead of \textquotedblleft\ $+$ \textquotedblright $,$
then 
\begin{equation}
\left\vert J\right\vert \leq \frac{1}{4}\left\vert Z-z\right\vert ^{2}.
\label{5.16}
\end{equation}
\end{corollary}

\begin{remark}
It is obvious that if one chooses the discrete measure above, then all the
inequalities in this section may be written for sequences of real or complex
numbers. We omit the details.
\end{remark}

\end{document}